\documentclass[a4paper,10pt]{article}
\usepackage{latexsym,amssymb,amscd}

\def\NZQ{\Bbb}               
\def\NN{{\NZQ N}}

\def\KK{{\NZQ K}}

\def\B'c{{\mathcal{B'}}}
\def\U'c{{\mathcal{U'}}}
\def\Mc{{\mathcal{M}}}

\def\Fc{{\mathcal{F}}}

\def\Vc{{\mathcal{V}}}

\def\ab{{\bold a}}

\def\kb{{\bold k}}

\def\poly#1#2#3{#1[#2_1,\dots,#2_{#3}]}

\newtheorem{Theorem}{Theorem}[section]

\newtheorem{Proposition}[Theorem]{Proposition}

\newtheorem{Remarks}[Theorem]{Remarks}

\newtheorem{Definition}[Theorem]{Definition}

\title{Computing Koszul Homology for Monomial Ideals}
\author{Eduardo S\'aenz de Cabez\'on\footnote {Partially supported by NEST Project 5006 (GIFT)}}

\begin{document}

\maketitle

\begin{abstract}
The Koszul homology of modules of the polynomial ring $R$ is a central object in commutative algebra. It is strongly related with the minimal free resolution of these modules, and thus with regularity, Hilbert functions, etc. Here we consider the case of modules of the form $R/I$ where $I$ is a monomial ideal. So far, some good algorithms have been given in the literature and implemented in different Computer Algebra Systems (e.g. CoCoa, Singular, Macaulay), which compute minimal free resolutions of modules of the form $R/I$ with $I$ an ideal in $R$, which include the case of $I$ being a monomial ideal as a particular one (a good review is given in \cite {Sie}). Our goal is to build algorithms specially targeted to monomial ideals, taking into account the particular combinatorial and structural properties of these ideals. This being a first goal, it is also a first step of an alternative approach to the computation of the Koszul homology and minimal free resolutions of polynomial ideals. 
\end{abstract}

\section{The Koszul Complex}
Let $\Vc$ be a vector space of dimension $n$ and give it the basis $\{x_1,\dots,x_n\}$ so that we can identify $SV$ with $R=\poly \kb x n$,  and consider the complex
$$\KK: 0\rightarrow R\otimes\wedge^n\Vc\stackrel{\partial}{\rightarrow}R\otimes\wedge^{n-1}\Vc\stackrel{\partial}{\rightarrow}\cdots R\otimes\wedge^{1}\Vc\stackrel{\partial}{\rightarrow}R\otimes\wedge^0\Vc\rightarrow\kb\rightarrow 0$$
where the maps are given by the rule
$$\partial(w_1\cdots w_q\otimes x_{i_0}\wedge\cdots\wedge x_{i_l})\longmapsto \sum_{j=0}^l (-1)^{j+1}w_1\cdots w_qx_{i_j}\otimes  x_{i_0}\wedge\cdots\wedge\widehat{x_{i_l}}\cdots\wedge x_{i_l}$$
this is called the {\it Koszul Complex}, and it is a minimal free resolution of $\kb$.

Given a graded module $\Mc$, its {\it Koszul Complex} $(\KK(\Mc), \partial)$ is the tensor product complex $\Mc\otimes_R \KK$. The {\it Koszul homology} of $\Mc$ is the homology of this product complex. We can identify this homology modules with $Tor_\bullet^R(\Mc,\kb)$, which can also be computed using any resolution of $M$ and tensoring with $\kb$.

\section{Combinatorial and Simplical nature of Monomial Ideals }
In the next pages, let $I\subseteq R=\poly \kb x n$ be a monomial ideal of the polynomial ring in $n$ variables over a field $\kb$. Let $\{m_1, \dots, m_r\}$ be the minimal generating set of $I$. For monomial ideals, Taylor \cite{Tay} gave a resolution that can be used for computing Koszul homology. This resolution has size $2^r$ and length $r$, which makes it too big for actual computations.  So, one of our first goals will be to reduce the size of the resolution, by computing it step by step and making use of the properties of the multigraded Betti numbers for monomial ideals.
\subsection{The LCM-Lattice and local computations of Koszul Homology}
 We denote (see \cite{GPW1}) by $L_I$ the lattice with elements labeled by the least common multiples of $m_1,\dots, m_r$ ordered by divisibility. We call $L_I$ the {\it LCM-Lattice} of $I$. We know, from Taylor's resolution that $\beta_{i,\ab}(I)=0$ if $\ab\notin L_I$.

On the other hand, if we build the minimal syzygy resolution step by step, we can assume that we have a minimal generating set of $Syz_i(I)$, and if we consider the multidegrees of these generators, the following proposition gives us a relationship between them and the multidegrees of the elements in a minimal generating set of $Syz_{i+1}(I)$:
\begin{Proposition}
Let $I\subset \poly \kb x n$ be a monomial ideal, and let $G_{i-1}$ be a minimal generating set of $Syz_{i_-1}(I)$. Let $md(G_{i-1})=\{ {\rm multideg}(s)\vert s\in G_{i-1}\}$ the set of multidegrees of the elements in $G_{i-1}$; then the set $md(G_{i})$ of multidegrees of the minimal generators of $Syz_i(I)$ is a subset of $\{{\rm lcm}(g_k,g_l)\vert g_k,g_l\in md(G_{i-1})\}$ 
\end{Proposition}

\subsection{Simplicial Koszul Complexes}
One way of computing the Koszul homology of $I$ at a given multidegree $\ab$ is to associate a simplicial complex to the ideal at that multidegree and express the Koszul homology of $I$ at $\ab$ in terms of this simplicial complex. See for example \cite{Bay}, \cite{P02} or \cite{MS04}. The basics are as follows:

Let $I$ be a monomial ideal and $x^{\ab}\in I$ let $l=|\mbox{support}(\ab)|$ i.e. $l$ is the number of nonzero coordinates of ${\ab}$. We associate to $\ab$ the $l$-simplex $\Delta_{\ab}$, where the vertices are labeled by the $a_i$ different from zero. Now, we build the subcomplex of $\Delta_{\ab}$ given by 
$$\Delta_{\ab}^I:=\{\tau \in \Delta_{\ab}\vert x^{\ab}/x^{ \tau}\in I\}$$ where $x^\tau$ is a squarefree monomial with exponents given by the variables defining the face $\tau$.

With this definition we have the following result (see \cite{Bay},\cite{MS04})
\begin{Proposition}
$$H_i(\KK(I)_\ab)\simeq \widetilde {H}_i(\Delta_{\ab}^I) \forall i$$
\end{Proposition}

\begin {Remarks}

\begin{enumerate}
\item{The correspondence of the chains of $S_\ab^I$ with chains of $\KK(I)_\ab$ that makes explicit this isomorphism  is just the $\kb$-linear extension of the correspondence that associates the face $\tau \in \Delta_{\ab}^I$ with $x^{\ab}/x^\tau\otimes x^{(\tau)}$ where, if $\tau=a_{\tau_1},\dots a_{\tau_k}$  the right-hand side of the latter has to be understood as $x^{(\tau)}=x^{a_{\tau_1}}\wedge\dots\wedge x^{ a_{\tau_k}}$.}
\item{Note that in the framework of our strategy for computing $H(\KK(I))_*$, in each step we are only interested on the computation of $H_i(\KK(I)_{\ab})$ for some ${\ab}$ and thus, we are only interested in one homological dimension at a time, so in particular we will only need a short piece of $\widetilde {C}_\bullet(\Delta_{\ab}^I;\kb)$ and we can use a subcomplex of $\Delta_{\ab}^I$, namely the one which has as facets the $i+1$ and $i$-faces of $\Delta_{\ab}^I$, since both homologies will be equal at dimension $i$.}
\end{enumerate}
\end{Remarks}

\section {Recursive Computations using Exact Sequences}
Given a monomial ideal $I\in R$ in minimally generated by $\{m_1,\dots,m_r\}$, let $I'$ be the ideal minimally generated by $m_1,\dots,m_{r-1}$ and $\tilde{I}$ the  ideal generated by $lcm(m_1,m_r),\dots,lcm(m_{r-1},m_r)$. With this notation:

\begin{Proposition}
The following sequence of complexes is exact:
$$0\longrightarrow \KK(\tilde{I}  )\longrightarrow \KK(I')\oplus \KK(\langle  m_r\rangle  )\longrightarrow \KK(I)\longrightarrow 0 $$
And as Koszul differential preserves total multidegree, we have the multigraded version, which is also exact:
$$0\longrightarrow \KK_{\ab}(\tilde{I}  )\longrightarrow \KK_{\ab}(I')\oplus \KK_{\ab}(\langle  m_r\rangle  )\longrightarrow \KK(I)_{\ab}\longrightarrow 0 $$ for all $\ab\in\NN^n$
\end {Proposition}

Now, given this short exact sequence of complexes, it gives raise to a long exact sequence in homology, by means of a connecting homomorphism $\Delta$:
$$\cdots\longrightarrow H_{i+1}(\KK_{\ab}(I))\stackrel{\Delta}{\longrightarrow}H_{i}(\KK_{\ab}(\tilde{I}  )\longrightarrow $$ $$H_i(\KK_{\ab}(I)\oplus \KK_{\ab}(\langle  m_r\rangle  ))\longrightarrow H_i(\KK_{\ab}(I))\stackrel{\Delta}{\longrightarrow}\cdots$$

Using recursively these exact sequences at every $\ab\in\NN^n$ we can compute the Koszul homology of $I$ from the homology of ideals with smaller minimal generating sets. In particular, the long exact sequence in homology located at the exact multidegrees in which we need actual computations, will be very frequently short sequences (from which we can automatically read the Betti numbers) or very short (i.e. only two nonzero elements) so that we have isomorphisms from which we can read the actual generators of the homology modules. This will indeed be the situation for some important families of monomial ideals (e.g. generic ideals). 

Similar although different considerations have been made for the polynomial case \cite{Sie} and for the monomial case \cite{Ser} as a strategy for the computations of minimal resolutions or Koszul homology, without combining it with the combinatorial or simplicial properties of monomial ideals. 

\section{Some special families of monomial ideals}
\subsection {Generic monomial ideals}
\begin{Definition}
We say that a monomial ideal $I=\langle m_1,\dots,m_r\rangle$ is {\rm generic} if whenever two minimal generators $m_i$ and $m_j$ have the same positive degree in some variable, a third generator $m_k$ strictly divides their least common multiple $lcm(m_i,m_j)$ (\cite{MS04} def.6.5).

Let $I$ be a monomial ideal with minimaly generating set $\{m_1,\dots,m_r\}$. The {\rm Scarf Complex} $\Delta_I$ is defined as
$$\Delta_I=\{\sigma\subseteq\{1,\dots,r\}\vert m_\sigma=m_\tau\Rightarrow\sigma=\tau\}$$
and the {\rm Algebraic Scarf Complex} of $I$ is the Taylor complex $\Fc_{\Delta_I}$ supported on $\Delta_I$ (cf.ibid. defs 6.7 and 6.11).
\end{Definition}
By means of the special relation between the Scarf Complex and generic Monomial Ideals (we know that $\Fc_{\Delta_I}$ minimally resolves $I$ if $I$ is generic), we have the following proposition:
\begin{Proposition}
If $I$ is a generic monomial ideal, then for every $\ab\in\NN$, the exact sequence $\cdots H_i(\KK(\tilde{I}))_\ab\rightarrow H_i(\KK(I'))_\ab\oplus H_i(\KK(\langle m_r\rangle))_\ab\rightarrow H_i(\KK(I))_\ab\rightarrow H_{i-1}(\KK(\tilde{I}))_\ab\rightarrow\cdots$
has at most two nonzero elements.
\end{Proposition}
\subsection {Quasi-stable monomial ideals}

We say that a monomial ideal $I$ is {\it quasi-stable} if it posseses a (finite) Pommaret basis \cite{Sei02a}. In \cite{Sei02b}, Seiler shows that provided we know a Pommaret basis of a module $\Mc$ of  $R=\poly \kb x n$, we obtain a closed form for a resolution of $\Mc$ of minimal length, which we call the {\it Pommaret-Seiler} resolution of $\Mc$ (P-S for short), and which is a generalization of the construction of Eliahou-Kervaire for stable ideals \cite{EK}. 
In the case of monomial ideals, the length of Taylor resolution is the number of generators given to describe the ideal $I$. This length is usually greater than {\it pdim}$(I)$, but normally, the size of the first modules in this resolution (we will denote it T) is much smaller than the corresponding size of the modules in P-S. Note that when we refer here to Taylor's resolution, we always consider the one generated by the minimal generating set of $I$.
The knowledge of the structure of a Pommaret basis of $I$ allows us to reduce P-S resolution so that the modules in the new resolution are smaller than the corresponding ones in T with minimal generating set. Thus, we obtain a subresolution of Taylor's which is of minimal length. This we can take as starting point of our considerations.

\section{Algorithm, examples and results}
An algorithm for computing Koszul Homology for monomial ideals can be constructed making use of the approaches we have seen in the precedent sections. First of all, we can reduce the number of corners (multidegrees) in which we will actually compute homology. This can be performed by combining the combinatorial properties of the $lcm$-lattice and the short exact sequences of Koszul complexes we have seen above. A second moment is the actual computation of the generators of the homology modules. Here we can make use of both the long exact sequence in Koszul homology we described, and the simplicial Koszul complex at a given multidegree.

The table below shows some features of the algorithm: We see how, applying the combinatorial properties of the monomial ideals, the number of multidegrees to be considered is rather low, compared with Taylor's resolution, and in some cases fairly close to the actual minimal resolution. The table shows the results for random monomial ideals $I$ in $n$ variables, with $g$ generators the degree of which varies between $Md$ (maximal degree of a generator of $I$) and $md$ (minimal degree of a generator of $I$):

\begin {table}[h]\caption{Number of multidegrees checked, compared with Taylor and minimal resolutions}
\begin {center}
    \begin{tabular}{|c|c|c||l|l|l|}
    \hline
    \small{{\bf n}}& \small{{\bf g}}& \small{{\bf Md/md}} & \small{{\bf Taylor}}
    &\small{{\bf Algorithm}}& \small{{\bf Minimal}}\\
    \hline
    3&6&35/18&64&18&17\\
    3&10&34/18&1024&36&35\\
    6&5&69/20&32&19&17\\
    6&12&80/25&4096&87&57\\
    6&15&82/39&32768&309&185\\
    9&7&112/68&128&109&81\\
    9&8&111/84&256&172&123\\
    \hline
    \end{tabular}

\end{center}
\end{table}

If we know that our ideal belongs to some special type, we can take advantage of this fact and obtain better results or even ad-hoc algorithms. This applies in particular to generic or quasi-stable ideals. But also some other strategies copuld be taken into account, for example, ideals with a big number of generators, compared with the number of variables, will tend to eliminate much more corners due to topological-simplicial reasons, and thus this strategy will be sometimes prefered to the recursion on the number of generators.

Next step is the integration of these monomial techniques into the context of computing Koszul Homology of polynomial ideals. Homological Perturbation Theory provides a good framework for this integration, that can be an alternative to the usual methods, see \cite{Ser},\cite{LS}.

\bigskip

\noindent
Eduardo S\'aenz de Cabez\'on\\
Departamento de Matematicas y Computaci\'on\\
Universidad de La Rioja\\
Edificio Vives, c/Luis de Ulloa s/n\\
26005 Logro\~no, La Rioja, Spain\\
e-mail:eduardo.saenz-de-cabezon@dmc.unirioja.es
\end{document}